\documentclass[12pt]{article}
\usepackage{amsmath,amssymb,amsthm,here,bm,framed,empheq}
\usepackage[all]{xy}
\usepackage[right=0.75in,left=0.75in,top=1in,bottom=1in]{geometry}

\newtheorem{theorem}{Theorem}[section]
\newtheorem{lemma}{Lemma}[section]
\newtheorem{proposition}{Proposition}[section]

\newtheorem{remark}{Remark}[section]
\newtheorem{example}{Example}[section]
\newtheorem{problem}{Problem}[section]

\DeclareMathOperator{\re}{{\rm Re}}

\newcommand{\abs}[1]{\left\vert{#1}\right\vert}

\newcommand{\GL}{{\rm GL}}

\makeatletter
  
  \@addtoreset{equation}{section}
  	
	\@addtoreset{figure}{section}

\title{\bf \vspace{-22pt}On a Connection Problem for the Generalized Hypergeometric Equation}
\date{\today}
\author{Shunya ADACHI}
\begin{document} 
\date{}
\vspace{-22pt}
\maketitle
\renewcommand{\thefootnote}{\fnsymbol{footnote}}
\renewcommand{\thefootnote}{\arabic{footnote}}
\vspace{-24pt}

\begin{abstract}
We study a connection problem between the fundamental systems of solutions at singular points $0$ and $1$ for the generalized hypergeometric equation which is satisfied by the generalized hypergeometric series ${}_nF_{n-1}$. In general, the local solution space around $x=1$ consists of one dimensional singular solution space and $n-1$ dimensional holomorphic solution space. Therefore in the case of $n\ge3$, the expression of connection matrix depends on the choice of the fundamental system of solutions at $x=1$.
On the connection problem for ordinary differential equations, Sch\"afke and Schmidt (LNM {\bf 810}, Springer, 1980) gave an impressive idea which focuses on the series expansion of fundamental system of solutions. We apply their idea to solve the connection problem for the generalized hypergeometric equation and derive the connection matrix. 

\end{abstract}

\section{Introduction}
We consider the {\it generalized hypergeometric differential equation}:
\begin{equation}\tag{E}\label{Eq}
\left[\delta\prod_{i=1}^{n-1}(\delta+\beta_i-1)-x\prod_{i=1}^n(\delta+\alpha_i)\right]y=0,
\end{equation}
where $x, \alpha_i,\beta_i\in\mathbb{C}$ and $\displaystyle\delta=x\frac{d}{dx}$. This equation has regular singular points at $x=0,1,\infty$ and whose rank is $n$. The Riemann scheme of \eqref{Eq} is given by
\begin{equation}\label{RS}
\left\{
\begin{array}{ccc}
x=0&x=1&x=\infty\\[1pt]
0&0&\alpha_1\\[1pt]
1-\beta_1&1&\alpha_2\\
\vdots&\vdots&\vdots\\
1-\beta_{n-2}&n-2&\alpha_{n-1}\\[1pt]
1-\beta_{n-1}&-\beta_n&\alpha_n
\end{array}
\right\}
\end{equation}
where $\beta_n$ is defined by $\alpha_1+\cdots+\alpha_n=\beta_1+\cdots+\beta_n$. The generalized hypergeometric equation is rigid, i.e., The linear ordinary differential equation whose Riemann scheme is given by \eqref{RS} is only the Generalized hypergeometric equation \eqref{Eq}. 

Throughout this paper, we assume 
\begin{equation}\label{nonresonance}
\beta_i,\,\beta_i-\beta_j\notin\mathbb{Z}
\end{equation}
for $1\le i,j\le n$ with $i\neq j$. Under this assumption, a fundamental system of solutions (i.e., basis of local solution space) around the origin is given by
\[\label{FS1}
\begin{aligned}
y_1^{[0]}(x)&={}_nF_{n-1}(\bm{\alpha}_0;\bm{\beta}_0;x),\\
y_{i+1}^{[0]}(x)&=x^{1-\beta_i}{}_nF_{n-1}(\bm{\alpha}_i;\bm{\beta}_i;x)
\end{aligned}
\]
where $\bm{\alpha}_0=(\alpha_1,\alpha_2,\ldots,\alpha_n),\,\bm{\beta}_0=(\beta_1,\ldots,\beta_{n-1})$ and 
\begin{align*}
&\bm{\alpha}_i=(\alpha_1+1-\beta_i,\alpha_2+1-\beta_i,\ldots,\alpha_n+1-\beta_i),\\
&\bm{\beta}_i=(\beta_1+1-\beta_i,\ldots,2-\beta_i,\ldots,\beta_{n-1}+1-\beta_i)
\end{align*}
for $1\le i\le n-1$. Here the symbol ${}_nF_{n-1}$ denotes the {\it generalized hypergeometric series}:
\begin{equation}\label{GHGS}
{}_nF_{n-1}\left(\begin{array}{c}
\alpha_1,\alpha_2,\ldots,\alpha_n\\
\beta_1,\ldots,\beta_{n-1}
\end{array};x\right)=\sum_{m\ge0}\frac{(\alpha_1)_m(\alpha_2)_m\cdots(\alpha_n)_m}{(\beta_1)_m\cdots(\beta_{n-1})_mm!}x^m.
\end{equation}
This series converges on $x\in D_0:=\{x\in\mathbb{C}\,;\,\abs{x}<1\}$ and the symbol $(a)_m$ is defined by
\[
(a)_m=\frac{\Gamma(a+m)}{\Gamma(a)}=\begin{cases}
1&m=0,\\
a(a+1)\cdots(a+m-1)&m\ge1.
\end{cases}
\]

Around $x=1$, a fundamental system of solutions is given by
\begin{equation}\label{FS2}
\begin{aligned}
y_{i}^{[1]}(x)&=(1-x)^{i-1}\sum_{m\ge0}d_m^{(i)}(1-x)^m,\quad d_0^{(i)}\neq0\\
y_{n}^{[1]}(x)&=(1-x)^{-\beta_n}\sum_{m\ge0}d_m^{(n)}(1-x)^m,\quad d_0^{(n)}\neq0
\end{aligned}
\end{equation}
for $1\le i\le n-1$. Here $y_1^{[1]}(x),y_2^{[1]}(x),\ldots,y_{n-1}^{[1]}(x)$ and the power series in $y_n^{[1]}$ are convergent on $D_1=\{x\in\mathbb{C}\,;\,\abs{1-x}<1\}$.
\begin{remark}\mbox{}\rm
\begin{enumerate}
\item To determine the fundamental system of solutions $y_{i}^{[1]}(x)\,(1\le i\le n-1)$, we have to fix the coefficients $\{d_{j}^{(i)}\}_{j=0}^{n-1-i}$. Then the other coefficients $\{d_{j}^{(i)}\}_{j\ge n-i}$ are determined uniquely. On the other hand, the solution $y_n^{[1]}(x)$ is determined by only fixing $d_n^{(0)}$.
\item In the case of $n=2$, the fundamental system of solutions \eqref{FS2} are expressed by Gauss hypergeometric series ${}_2F_1$. 
\end{enumerate}
\end{remark}

In this paper, we consider a connection problem for the fundamental systems of solutions at singular points $x=0$ and $1$. That is, we consider the following problem.
\begin{problem}[Connection problem]\label{Prob_1}
We set $n$-row vectors
\[
\mathcal{Y}^{[0]}=(y_1^{[0]},y_2^{[0]},\ldots,y_n^{[0]}),\qquad \mathcal{Y}^{[1]}=(y_1^{[1]},y_2^{[1]},\ldots,y_n^{[1]}).
\]
Then determine the connection matrix $C\in\GL(n,\mathbb{C})$ such that
\begin{equation}\label{CM}
\mathcal{Y}^{[0]}=\mathcal{Y}^{[1]}C,\quad x\in D_0\cap D_1\setminus\{0,1\}
\end{equation}
with $\arg{x}=\arg(1-x)=0$ on $0<x<1$.
\end{problem}

\medskip
Connection problem of the generalized hypergeometric equation between $x=0$ and $\infty$ is studied well (cf. Barnes \cite{Barnes1908}, Kawahata \cite{Kawahata1978}, Mimachi \cite{Mimachi2008, Mimachi2011}, Smith \cite{Smith1938}). Recently, {Oshima} \cite{Oshima2012} studied the connection problem for more general Fuchsian equations from the viewpoint of the theory of {middle convolution}.

On the connection problem of the generalized hypergeometric equation between $x=0$ and $1$, there is a difficulty that did not appear in the case of between $x=0$ and $\infty$: As is seen in \eqref{FS2}, the solution space around $x=1$ of \eqref{Eq} consists of one dimensional singular solution space and $n-1$ dimensional holomorphic solution space. Therefore in the case of $n\ge3$, the expression of connection matrix depends on the choice of the fundamental system of solutions at $x=1$. 

Mimachi \cite{Mimachi2008} considered the connection problem between singular solution at $x=1$ and fundamental system of solutions at $x=0$ in the case of $n=3$, and later he solved the same problem with the general case, in \cite{Mimachi2011}. Matsuhira-Nagoya \cite{MN2019} considered the connection problem between the fundamental systems of solutions and derived the connection matrix. 

In these studies, the connection problem was solved by using the integral representation of solutions and the notion of twisted cycle. That is, Matsuhira-Nagoya's result gives an answer to Problem \ref{Prob_1} in terms of the integral representation of solutions. 

In this paper, we solve Problem \ref{Prob_1} from the viewpoints of series representation of local solutions. In doing this, Sch\"{a}fke-Schmidt's work \cite{SS1980} plays an important role: They focused on the series representation of local solutions, and gave an impressive idea for connection problem for general linear differential equations.
\begin{theorem}[\cite{SS1980}] \label{Lem_SS}
We consider the linear differential equation
\[
p_0(x)y^{(n)}+p_1(x)y^{(n-1)}+\cdots+p_{n}(x)y=0
\]
which has regular singular points $x=0,1$ and doesn't have any other singular points on $\abs{x}\le 1$. If this equation has local solutions 
\[
y^{[0]}(x)=x^\alpha\sum_{m\ge0}a_mx^m
\]
on $D_0\setminus\{0\}$ and 
\[
y_j^{[1]}(x)=(1-x)^{\alpha_j}\sum_{m\ge0}d_m^{(j)}(1-x)^m,\quad1\le j\le n,
\]
on $D_1\setminus\{1\}$,
then the following asymptotic formula for $a_m$ is hold:
\begin{align*}
a_m&=\sum_{j=1}^n\frac{\Gamma(m+\alpha-\alpha_j)}{\Gamma(m+\alpha+1)}\left\{\sum_{\ell=0}^k\left(\prod_{s=1}^\ell\frac{-s-\alpha_j}{m+\alpha-s-\alpha_j}\right)d_\ell^{(j)}\right\}\frac{c_j}{\Gamma(-\alpha_j)}+O(m^{-\alpha_--k-2}), \quad m\to\infty
\end{align*}
where $c_1,c_2,\ldots,c_n\in\mathbb{C}$ are the connection coefficients such that
\[
y^{[0]}(x)=\sum_{j=1}^nc_jy_j^{[1]}(x),\quad 0<x<1
\]
with $\arg{x}=\arg(1-x)=0$, $k\in\mathbb{N}$ is arbitrary number and
\[
\alpha_-=\min\,\{\re\alpha_j\,;\,j~{\rm s.t.}~\alpha_j\notin\mathbb{N}\}.
\]
\end{theorem}

\begin{remark}\rm In \cite{SS1980}, the authors considered linear systems of first order ordinary differential equations. But the above statement can be shown by a similar way.
\end{remark}
Theorem \ref{Lem_SS} says that the connection coefficients appear in the asymptotic behavior of the coefficients of local solutions. We use this theorem to solve the connection problem for the generalized hypergeometric equation and derive the connection matrix explicitly. Our main result gives an example which shows Sch\"{a}fke and Schmidt's result is useful for derivation of the connection matrix of the fundamental systems of solutions for Fuchsian differential equations.

\medskip

This paper is organized as follows. In Section \ref{section_main} we give our main result and the strategy of a proof. Since Theorem \ref{main} is obtained from Proposition \ref{main2}, we show Proposition \ref{main2} in Section \ref{section_proof3.2}. Section \ref{section_proofprop} is devoted to proofs of Propositions \ref{c1-cn-1} and \ref{hgs_conv} which are necessary for our proof of Proposition \ref{main2} are given. 

\section{Main Result}\label{section_main}
In this section, we state our main theorem. Let us set
\[
D=
\left(\begin{array}{cccccc}d_0^{(1)} &  &  &  &  &  \\[3pt]d_1^{(1)} & d_0^{(2)} &  &  &  &  \\[3pt]d_2^{(1)} & d_1^{(2)} & d_0^{(3)} &  &  &  \\[3pt]\vdots & \vdots & \vdots & \ddots &  &  \\[3pt]d_{n-2}^{(1)} & d_{n-3}^{(2)} & d_{n-4}^{(3)} & \cdots & d_0^{(n-1)} &  \\[3pt]0 & 0 & 0 & \cdots & 0 & d_0^{(n)}\end{array}\right).
\]
Since we assumed $d_0^{(1)}d_0^{(2)}\cdots d_0^{(n)}\neq0$, the matrix $D$ is invertible. Our main result is as follows.
\begin{theorem}\label{main}
Assume that the condition \eqref{nonresonance} and $\re\beta_n<-n+2$ hold. Then the connection matrix $C$ in \eqref{CM} is given by
\[
C=D^{-1}P,
\]
where $P=(p_{ij})_{1\le i,j\le n}$ is the $n\times n$ matrix whose entries are
\[
p_{ij}=\begin{cases}
\dfrac{(-1)^{i-1}}{(i-1)!}\dfrac{(\bm{\alpha}_{j-1})_{i-1}}{(\bm{\beta}_{j-1})_{i-1}}{}_nF_{n-1}\left(
\bm{\alpha}_{j-1}+i-1;\bm{\beta}_{j-1}+i-1;1\right)
&1\le i\le n-1,\\[11pt]
\dfrac{\Gamma(\bm{\beta}_{j-1},\beta_n)}{\Gamma(\bm{\alpha}_{j-1})}&i=n.
\end{cases}
\]
The symbols $(\bm{\alpha})_m$ and $\Gamma(\bm{\alpha})$ are defined by
\[
(\bm{\alpha})_m:=\prod_{i=1}^n(\alpha_i)_m,\quad \Gamma(\bm{\alpha}):=\prod_{i=1}^n\Gamma(\alpha_i)
\]
respectively, and 
\[
\bm{\alpha}+a:=(\alpha_1+a,\alpha_2+a,\ldots,\alpha_n+a)
\]
for $m\in\mathbb{N},\,a\in\mathbb{C}$ and $\bm{\alpha}=(\alpha_1,\alpha_2,\ldots,\alpha_n)\in\mathbb{C}^n$. 
\end{theorem}
\begin{remark}\rm
We note that this theorem holds even in the cases
\begin{enumerate}
\item $\alpha_i\in\mathbb{Z}$ for some $1\le i\le n$,
\item $\alpha_i-\beta_j\in\mathbb{Z}$ for some $1\le i\le n$ and $1\le j\le n-1$,
\item $\alpha_i-\alpha_j\in\mathbb{Z}$ for some $1\le i,j\le n$.
\end{enumerate} 
These cases were not treated in Matsuhira-Nagoya \cite{MN2019}, and it is known that in the cases (1) and (2), our equation \eqref{Eq} is reducible (cf. Beukers-Heckman \cite{BH1989}).  
\end{remark}

\begin{example}
In the case of $n=3$, 
\[
D=\begin{pmatrix}
d_0^{(1)}&&\\[2pt]
d_1^{(1)}&d_0^{(2)}&\\
0&0&d_0^{(3)}
\end{pmatrix}
\]
and{\small
\[
C=D^{-1}
\begin{pmatrix}
{}_3F_2(\bm{\alpha}_0;\bm{\beta}_0;1)&{}_3F_2(\bm{\alpha}_1;\bm{\beta}_1;1)&{}_3F_2(\bm{\alpha}_2;\bm{\beta}_2;1)\\[4pt]
-\dfrac{(\bm{\alpha}_0)_1}{(\bm{\beta}_0)_1}{}_3F_2(\bm{\alpha}_0+1;\bm{\beta}_0+1;1)&-\dfrac{(\bm{\alpha}_1)_1}{(\bm{\beta}_1)_1}{}_3F_2(\bm{\alpha}_1+1;\bm{\beta}_1+1;1)&-\dfrac{(\bm{\alpha}_2)_1}{(\bm{\beta}_2)_1}{}_3F_2(\bm{\alpha}_2+1;\bm{\beta}_2+1;1)\\[4pt]
\dfrac{\Gamma(\bm{\beta}_0,\beta_3)}{\Gamma(\bm{\alpha}_0)}&\dfrac{\Gamma(\bm{\beta}_1,\beta_3)}{\Gamma(\bm{\alpha}_1)}
&\dfrac{\Gamma(\bm{\beta}_2,\beta_3)}{\Gamma(\bm{\alpha}_2)}
\end{pmatrix}.
\]}
\begin{remark}\rm
In this case, the connection matrix is complicated because $D$ is not a diagonal matrix in general. But If we take fundamental system of solutions \eqref{FS2} satisfies $d_1^{(1)}=0$, or retaking $y_2(x)$ as $y_2(x)-d_1^1/d_0^1(1-x)y_1(x)$, we can make the matrix $D$ into diagonal matrix and the expression of connection matrix simple. 
\end{remark}
\end{example}

To give a proof of Theorem \ref{main}, it is sufficient to show the following proposition.
\begin{proposition}\label{main2}
Assume that the condition \eqref{nonresonance} and $\re\beta_n<-n+2$ hold. Then the connection coefficients $c_1,c_2,\ldots,c_n\in\mathbb{C}$ such that
\begin{equation}\label{eq_prop1_1}
y_1^{[0]}(x)=\sum_{j=1}^n c_jy_j^{[1]}(x),\quad x\in D_0\cap D_1\setminus\{0,1\}
\end{equation}
with $\arg x=\arg(1-x)=0$ on $0<x<1$ are given by
\begin{equation}\label{Dc}
\begin{pmatrix}
c_1\\
c_2\\
\vdots\\
c_{n-1}\\
c_n
\end{pmatrix}=
D^{-1}\begin{pmatrix}
{}_nF_{n-1}(\bm{\alpha}_0;\bm{\beta}_0;1)\\[2pt]
-\dfrac{(\bm{\alpha}_0)_{1}}{(\bm{\beta}_0)_{1}}{}_nF_{n-1}(\bm{\alpha}_0+1;\bm{\beta}_0+1;1)\\
\vdots\\
\dfrac{(-1)^{n-2}}{(n-2)!}\dfrac{(\bm{\alpha}_0)_{n-2}}{(\bm{\beta}_0)_{n-2}}{}_nF_{n-1}(\bm{\alpha}_0+n-2;\bm{\beta}_0+n-2;1)\\
\dfrac{\Gamma(\bm{\beta}_0,\beta_n)}{\Gamma(\bm{\alpha}_0)}
\end{pmatrix}.
\end{equation}
\end{proposition}
The expression \eqref{Dc} gives the first column of the connection matrix $C$ and other columns are obtained from the first column. This is seen from the fact that the Riemann scheme \eqref{RS} determines the generalized hypergeometric equation \eqref{Eq}. Let us explain this for detail. To clarify the dependence for parameters, we rewrite \eqref{eq_prop1_1} as
\[\label{depend}
y_{1}^{[0]}(x)=y_{1}^{[0]}(\bm{\alpha}_0;\bm{\beta}_0;x),\quad y_{j}^{[1]}(x)=y_{j}^{[1]}(\bm{\alpha}_0;\bm{\beta}_0;x),\quad  c_j=c_j(\bm{\alpha}_0;\bm{\beta}_0),\quad 1\le j\le n
\]
for a while. We note that the other solutions at $x=0$ is written as
\begin{equation*}\label{depend2}
y_{i+1}^{[0]}(\bm{\alpha}_0;\bm{\beta}_0;x)=x^{1-\beta_i}y_{1}^{[0]}(\bm{\alpha}_i;\bm{\beta}_i;x)
\end{equation*}
for $1\le i\le n-1$. In the following, we focus on the fundamental system of solutions at $x=1$ and show
\begin{equation}\label{tf1}
x^{1-\beta_i}y_{j}^{[1]}(\bm{\alpha}_i;\bm{\beta}_i;x)=y_{j}^{[1]}(\bm{\alpha}_0;\bm{\beta}_0;x),\quad 1\le j\le n.
\end{equation}
 Let us exchange the parameter
\begin{equation}\label{proc1}
(\bm{\alpha}_0,\bm{\beta}_0)\longmapsto (\bm{\alpha}_i,\bm{\beta}_i)
\end{equation}
in our equation \eqref{Eq}. Then we have the equation which is satisfied by the functions $y_{1}^{[0]}(\bm{\alpha}_i;\bm{\beta}_i;x)$ and $y_{j}^{[1]}(\bm{\alpha}_i;\bm{\beta}_i;x)$ and its Riemann scheme is given by
\begin{equation}\label{RS2}
\left\{
\begin{array}{ccc}
x=0&x=1&x=\infty\\[1pt]
0&0&\alpha_1+1-\beta_i\\[1pt]
\beta_i-\beta_1&1&\alpha_2+1-\beta_i\\
\vdots&\vdots&\vdots\\
\beta_i-1&i&\alpha_{i+1}+1-\beta_i\\
\vdots&\vdots&\vdots\\
\beta_i-\beta_{n-2}&n-2&\alpha_{n-1}+1-\beta_i\\[1pt]
\beta_i-\beta_{n-1}&-\beta_n&\alpha_n+1-\beta_i
\end{array}
\right\}.
\end{equation}
Next we consider the transform of unknown function
\begin{equation}\label{proc2}
y(x)\longmapsto x^{1-\beta_i}y(x).
\end{equation}
Then the Riemann scheme \eqref{RS2} turns into
\[\label{RS3}
\left\{
\begin{array}{ccc}
x=0&x=1&x=\infty\\[1pt]
1-\beta_i&0&\alpha_1\\[1pt]
1-\beta_1&1&\alpha_2\\
\vdots&\vdots&\vdots\\
0&i&\alpha_{i+1}\\
\vdots&\vdots&\vdots\\
1-\beta_{n-2}&n-2&\alpha_{n-1}\\[1pt]
1-\beta_{n-1}&-\beta_n&\alpha_n
\end{array}
\right\}.
\]
The equation which corresponds this scheme is nothing but \eqref{Eq}. This means the relation \eqref{tf1} holds. Therefore, we see that the connection formula \eqref{eq_prop1_1} is transformed into
\[
y_{i+1}^{[0]}(x)=c_1(\bm{\alpha}_i;\bm{\beta}_i)y_{1}^{[1]}(\bm{\alpha}_0;\bm{\beta}_0;x)+c_2(\bm{\alpha}_i;\bm{\beta}_i)y_{2}^{[1]}(\bm{\alpha}_0;\bm{\beta}_0;x)+\cdots+c_n(\bm{\alpha}_i;\bm{\beta}_i)y_{n}^{[1]}(\bm{\alpha}_0;\bm{\beta}_0;x)
\]
by following the procedures \eqref{proc1} and \eqref{proc2}. This gives the $i$-th column of the connection matrix $C$.

\section{Proof of Proposition \ref{main2}}\label{section_proof3.2}
We shall show Proposition \ref{main2} by dividing the proof into the following two steps.
\begin{description}
\item[\bf Step 1.] To derive 
\begin{equation}\label{cn}
d_0^{(n)}c_n=\dfrac{\Gamma(\bm{\beta}_0,\beta_n)}{\Gamma(\bm{\alpha}_0)}.
\end{equation}
\item[\bf Step 2.] To derive 
\begin{equation}\label{rel_1}
\sum_{j=1}^id_{i-j}^{(j)}c_j=
\frac{(-1)^{i-1}}{(i-1)!}\frac{(\bm{\alpha}_0)_{i-1}}{(\bm{\beta}_0)_{i-1}}{}_nF_{n-1}(\bm{\alpha}_0+i-1,\bm{\beta}_0+i-1;1)
\end{equation}
for $1\le i\le n-1$.
\end{description}
The relation \eqref{rel_1} gives the relation of $c_1,c_2,\ldots,c_{n-1}$ in \eqref{Dc}: The relation \eqref{rel_1} is rewritten as 
\begin{align*}
d_0^{(1)}c_1&={}_nF_{n-1}(\bm{\alpha}_0,\bm{\beta}_0;1),\\
d_0^{(2)}c_2+d_1^{(1)}c_1&=-\frac{(\bm{\alpha}_0)_{1}}{(\bm{\beta}_0)_{1}}{}_nF_{n-1}(\bm{\alpha}_0+1,\bm{\beta}_0+1;1),\\
&~\vdots\\
d_0^{(n-1)}c_{n-1}+\cdots+d_{n-2}^{(1)}c_1&=\frac{(-1)^{n-2}}{(n-2)!}\frac{(\bm{\alpha}_0)_{n-2}}{(\bm{\beta}_0)_{n-2}}{}_nF_{n-1}(\bm{\alpha}_0+n-2;\bm{\beta}_0+n-2;1).
\end{align*}
By expressing them as a matrix form, we have
\[
\left(\begin{array}{ccccc}d_0^{(1)} &  &   &  \\[3pt]
d_1^{(1)} & d_0^{(2)} &  &    \\[3pt]
\vdots & \vdots & \ddots &     \\[3pt]
d_{n-2}^{(1)} & d_{n-3}^{(2)} &  \cdots & d_0^{(n-1)} 
\end{array}\right)
\begin{pmatrix}
c_1\\
c_2\\
\vdots\\
c_{n-1}
\end{pmatrix}=
\begin{pmatrix}
{}_nF_{n-1}(\bm{\alpha}_0;\bm{\beta}_0;1)\\[2pt]
-\dfrac{(\bm{\alpha}_0)_{1}}{(\bm{\beta}_0)_{1}}{}_nF_{n-1}(\bm{\alpha}_0+1;\bm{\beta}_0+1;1)\\
\vdots\\
\dfrac{(-1)^{n-2}}{(n-2)!}\dfrac{(\bm{\alpha}_0)_{n-2}}{(\bm{\beta}_0)_{n-2}}{}_nF_{n-1}(\bm{\alpha}_0+n-2;\bm{\beta}_0+n-2;1)
\end{pmatrix}.
\]

\subsection{Step 1 -Derivation of \eqref{cn}-}
In this subsection, we derive \eqref{cn}.
 By setting
\[
y_{1}^{[0]}(x)={}_nF_{n-1}(\bm{\alpha}_0;\bm{\beta}_0;x)=\sum_{m\ge0}a_mx^m,\quad a_m=\frac{(\alpha_1)_m(\alpha_2)_m\cdots(\alpha_n)_m}{(\beta_1)_m(\beta_2)_m\cdots(\beta_{n-1})m!}
\]
and applying Theorem \ref{Lem_SS}, we have
\[a_m=\frac{\Gamma(m+\beta_n)}{\Gamma(m+1)}\left\{\sum_{\ell=0}^k\left(\prod_{s=1}^\ell\frac{-s+\beta_n}{m-s+\beta_n}\right)d_\ell^{(n)}\right\}\frac{c_n}{\Gamma(\beta_n)}+O(m^{\re\beta_n-k-2}),\quad m\to\infty.
\]
Taking $k=0$ in this expression, we get
\[
a_m=\frac{\Gamma(m+\beta_n)}{\Gamma(m+1)}d_0^n\frac{c_n}{\Gamma(\beta_n)}+O(m^{\re\beta_n-2}).
\]
From this and the formula 
\[
\frac{\Gamma(m+1)}{\Gamma(m+\beta_n)}=m^{1-\beta_n}(1+O(m^{-1}))
\]
which is obtained by Stirling's formula, we obtain
\begin{equation}\label{cn_1}
d_0^{(n)}c_n=\frac{\Gamma(\beta_n)\Gamma(m+1)}{\Gamma(m+\beta_n)}a_m+O(m^{-1}).
\end{equation}
By substituting 
\begin{align*}
a_m&=\frac{(\alpha_1)_m(\alpha_2)_m\cdots(\alpha_n)_m}{(\beta_1)_m(\beta_2)_m\cdots(\beta_{n-1})m!}\\
&=\frac{\Gamma(\alpha_1+m)\Gamma(\alpha_2+m)\cdots\Gamma(\alpha_n+m)}{\Gamma(\beta_1+m)\Gamma(\beta_2+m)\cdots\Gamma(\beta_{n-1}+m)\Gamma(1+m)}\frac{\Gamma(\beta_1)\Gamma(\beta_2)\cdots\Gamma(\beta_{n-1})}{\Gamma(\alpha_1)\cdots\Gamma(\alpha_n)}
\end{align*}
into \eqref{cn_1} we have
\begin{equation}\label{cn_2}
d_0^{(n)}c_n=\frac{\Gamma(\bm{\beta}_0,\beta_n)}{\Gamma(\bm{\alpha}_0)}\frac{\Gamma(\alpha_1+m,\alpha_2+m,\ldots,\alpha_n+m)}{\Gamma(\beta_1+m,\ldots,\beta_{n-1}+m,\beta_n+m)}+O(m^{-1}).
\end{equation}
From Stirling's formula it holds that
\[
\frac{\Gamma(\alpha_1+m,\alpha_2+m,\ldots,\alpha_n+m)}{\Gamma(\beta_1+m,\ldots,\beta_{n-1}+m,\beta_n+m)}\to m^{\alpha_1+\alpha_2+\cdots+\alpha_n-(\beta_1+\cdots+\beta_{n-1}+{\beta_n})}=1
\]
as $m\to\infty$. Therefore by taking a limit of $m\to\infty$ for \eqref{cn_2}, we get the conclution\eqref{cn}.

\subsection{Step 2 -Derivation of \eqref{rel_1}-}
The relation \eqref{rel_1} is obtained from the following two propositions.
\begin{proposition}\label{c1-cn-1}
It holds that 
\begin{equation}\label{eq_c1-cn-1}
\sum_{j=1}^id_{i-j}^{(j)}c_j=
\frac{(-1)^{i-1}}{(i-1)!}\sum_{h=0}^m[h]_{i-1}a_h+O(m^{\re\beta_n+i-1}),\quad m\to\infty
\end{equation}
for $1\le i\le n-1$. Here we set
\[
[h]_{j}=\begin{cases}
1&j=0,\\
h(h-1)(h-2)\cdots(h-j+1)&j\ge1.
\end{cases}
\]
\end{proposition}
\begin{proposition}\label{hgs_conv}
Assume $\re\beta_n<-n+2$. Then we have
\begin{equation}\label{eq_hgs_conv}
\lim_{m\to\infty}\sum_{h=0}^m[h]_{i-1}a_h=\frac{(\bm{\alpha}_0)_{i-1}}{(\bm{\beta}_0)_{i-1}}{}_nF_{n-1}(\bm{\alpha}_0+i-1,\bm{\beta}_0+i-1;1)
\end{equation}
for $1\le i\le n-1$.
\end{proposition}
Proofs of these propositions are given in the following section. We admit these propositions and finish the proof of Proposition \ref{main2}. From the assumption, for any $1\le i\le n-1$ it is obtained
\[
\re\beta_n+i-1<-n+2+i-1\le -n+2+(n-1)-1=0.
\] 
Therefore we see that $O(m^{\re\beta_n+i-1})\to0$ as $m\to\infty$ for each $1\le i\le n-1$. Then by taking a limit $m\to\infty$ in \eqref{eq_c1-cn-1} and using \eqref{eq_hgs_conv}, we get the desired expression \eqref{rel_1}.

\section{Proofs of Propositions \ref{c1-cn-1} and \ref{hgs_conv}}\label{section_proofprop}
\subsection{Proof of Proposition \ref{c1-cn-1}}
We proof Proposition \ref{c1-cn-1} by induction on $i$. At first, we show the case of $i=1$, i.e.,
\begin{equation}\label{c1}
d_0^{(1)}c_1=\sum_{h=0}^{m}a_h+O(m^{\re\beta_n}),\quad m\to\infty.
\end{equation}
Let us consider the gauge transform 
\[
z(x)=(1-x)^{-1}y(x)
\]
in our equation \eqref{Eq}. Then the equation which is satisfied by $z(x)$ is given by 
\begin{equation}\label{RS4}
\left\{
\begin{array}{ccc}
x=0&x=1&x=\infty\\[1pt]
0&-1&\alpha_1+1\\[1pt]
1-\beta_1&0&\alpha_2+1\\
\vdots&\vdots&\vdots\\
1-\beta_{n-2}&n-3&\alpha_{n-1}+1\\[1pt]
1-\beta_{n-1}&-\beta_n-1&\alpha_n+1
\end{array}
\right\}.
\end{equation}
The fundamental system of solutions of this equation is obtained by considering the gauge transform for each solution of the equation \eqref{Eq} (given in \eqref{FS2}). We write them as
\begin{align*}
z^{[0]}(x)&=(1-x)^{-1}y_1^{[0]}(x)=\sum_{m\ge0}\left(\sum_{h=0}^ma_h\right)x^m,\\
z_j^{[1]}(x)&=(1-x)^{-1}y_j^{[1]}(x)=(1-x)^{j-2}\sum_{m\ge0} d_m^{(j)}(1-x)^m,\quad 1\le j \le n-1,\\
z_n^{[1]}(x)&=(1-x)^{-1}y_n^{[1]}(x)=(1-x)^{-\beta_n-1}\sum_{m\ge0} d_m^{(n)}(1-x)^m.
\end{align*}
Then the connection formula \eqref{eq_prop1_1} turns into
\[
z^{[0]}(x)=c_1z_1^{[1]}(x)+c_2z_2^{[1]}(x)+\cdots+c_nz_n^{[1]}(x).
\]
Here we remark that $z^{[0]}(x)$ is a holomorphic solution of \eqref{RS4} and its characteristic exponent is 0. The function $z_1^{[1]}(x)$ is a singular solution and its characteristic exponent is $-1$.

Now we use Theorem \ref{Lem_SS} for the equation \eqref{RS4}. Then we have
\begin{align*}
\sum_{h=0}^{m}a_h&=\left\{\sum_{\ell=0}^k\left(\prod_{s=1}^\ell\frac{-s+1}{m-s+1}\right)d_\ell^{(1)}\right\}\frac{c_1}{\Gamma(1)}\\
&+\frac{\Gamma(m+\beta_n+1)}{\Gamma(m+1)}\left\{\sum_{\ell=0}^k\left(\prod_{s=1}^\ell\frac{-s+\beta_n+1}{m-s+\beta_n+1}\right)d_\ell^{(n)}\right\}\frac{c_n}{\Gamma(\beta_n+1)}+O(m^{-\alpha_--k-2})
\end{align*}
for any $k\in\mathbb{N}$. Since 
\[
\prod_{s=1}^\ell\frac{-s+1}{m-s+1}=0
\]
for all $\ell\ge1$, it holds that 
\[
\sum_{\ell=0}^k\left(\prod_{s=1}^\ell\frac{-s+1}{m-s+1}\right)d_\ell^{(1)}=\sum_{\ell=0}^0\left(\prod_{s=1}^\ell\frac{-s+1}{m-s+1}\right)d_\ell^{(1)}=d_0^{(1)}.
\]
Therefore we have
\begin{align*}
\sum_{h=0}^{m}a_h=d_0^{(1)}c_1+\frac{\Gamma(m+\beta_n+1)}{\Gamma(m+1)}\left\{\sum_{\ell=0}^k\left(\prod_{s=1}^\ell\frac{-s+\beta_n+1}{m-s+\beta_n+1}\right)d_\ell^{(n)}\right\}\frac{c_n}{\Gamma(\beta_n+1)}+O(m^{-\alpha_--k-2}).
\end{align*}
Here $\alpha_-=\min\{-1,-\re\beta_n-1\}$. Since we assumed $\re\beta_n<-n+2$ and $n\ge2$, we have
\[
-\re\beta_n-1>n-3\ge -1.
\]
Therefore we see that $\alpha_-=-1$. Now we fix $k\in\mathbb{N}$ which satisfies $k\ge -\re\beta_n-1$. Then it holds
\[
-\alpha_--k-2=-k-1\le \re\beta_n.
\]
From this we have
\[
O(m^{-\alpha_--k-2})=O(m^{-k-1})=O(m^{\re\beta_n}).
\]
Next we focus on the second term of the right-hand side. From Stirling's formula, it holds that
\[
\frac{\Gamma(m+\beta_n+1)}{\Gamma(m+1)}=O(m^{\re\beta_n}).
\]
Therefore we have the desired relation \eqref{c1}.

Next let $i\in\mathbb{N}$ with $2\le i\le n-1$, and suppose that \eqref{eq_c1-cn-1} (with $i$ replaced by $p$) is already proved for all $p\in\mathbb{N}$ with $1\le p\le i-1$. We consider the gauge transform 
\[
z(x)=(1-x)^{-i}z(x)
\]
in the equation \eqref{Eq}.  Then the equation which is satisfied by $z(x)$ is given by 
\begin{equation}\label{RS5}
\left\{
\begin{array}{ccc}
x=0&x=1&x=\infty\\[1pt]
0&-i&\alpha_1+i\\[1pt]
1-\beta_1&-i+1&\alpha_2+i\\
\vdots&\vdots&\vdots\\
1-\beta_{i-1}&-1&\alpha_{i}+i \\
1-\beta_i&0&\alpha_{i+1}+i\\
\vdots&\vdots&\vdots\\
1-\beta_{n-1}&-\beta_n-i&\alpha_n+i
\end{array}
\right\}.
\end{equation} 
The fundamental system of solutions of this equation is obtained by considering the gauge transform for each solution of the equation \eqref{Eq} (given in \eqref{FS2}). We write them as
\begin{align*}
z^{[0]}(x)&=(1-x)^{-i}y_1^{[0]}(x)=\sum_{m\ge0}\left(\frac{1}{(i-1)!}\sum_{h=0}^m(m-h+1)_{i-1}a_h\right)x^m,\\
z_j^{[1]}(x)&=(1-x)^{-i}y_j^{[1]}(x)=(1-x)^{j-i-1}\sum_{m\ge0} d_m^{(j)}(1-x)^m,\quad 1\le j \le n-1,\\
z_n^{[1]}(x)&=(1-x)^{-i}y_n^{[1]}(x)=(1-x)^{-\beta_n-i}\sum_{m\ge0} d_m^{(n)}(1-x)^m.
\end{align*}
Then the connection formula \eqref{eq_prop1_1} turns into
\[
z^{[0]}(x)=c_1z_1^{[1]}(x)+c_2z_2^{[1]}(x)+\cdots+c_nz_n^{[1]}(x).
\]
Here we remark that $z^{[0]}(x)$ is a holomorphic solution of \eqref{RS5} and its characteristic exponent is 0. The functions $z_j^{[1]}(x)\,(1\le j\le i)$ are singular solutions of \eqref{RS5} and characteristic exponent of each solution is $j-i-1$.

Now we use Theorem \ref{Lem_SS} for the equation \eqref{RS5}. Then we have
\begin{align*}
&\frac{1}{(i-1)!}\sum_{h=0}^m(m-h+1)_{i-1}a_h\\
&=\sum_{j=1}^{i}\frac{\Gamma(m+i-j+1)}{\Gamma(m+1)}\left\{\sum_{\ell=0}^k\left(\prod_{s=1}^\ell\frac{-s+i-j+1}{m-s+i-j+1}\right)d_\ell^{(j)}\right\}\frac{c_j}{\Gamma(i-j+1)}\\
&\quad +\frac{\Gamma(\beta_n+i+m)}{\Gamma(m+1)}\left\{\sum_{\ell=0}^k\left(\prod_{s=1}^\ell\frac{-s+\beta_n+i}{m-s+\beta_n+i}\right)d_\ell^{(n)}\right\}\frac{c_n}{\Gamma(\beta_n+i)}+O(m^{-\alpha_--k-2})
\end{align*}
for any $k\in\mathbb{N}$. Now we fix $k\in\mathbb{N}$ which satisfies $k\ge i-1\,(\ge i-j+1)$. Since 
\[
\prod_{s=1}^\ell\frac{-s+i-j+1}{m-s+i-j+1}=0
\]
for all $\ell\ge i-j+1$, it holds that
\[
\sum_{\ell=0}^k\left(\prod_{s=1}^\ell\frac{-s+i-j+1}{m-s+i-j+1}\right)d_\ell^{(j)}=\sum_{\ell=0}^{i-j}\left(\prod_{s=1}^\ell\frac{-s+i-j+1}{m-s+i-j+1}\right)d_\ell^{(j)}.
\]
Therefore we have
\begin{equation}
\begin{aligned}\label{ind1}
&\frac{1}{(i-1)!}\sum_{h=0}^m(m-h+1)_{i-1}a_h\\
&=\sum_{j=1}^{i}\frac{\Gamma(m+i-j+1)}{\Gamma(m+1)}\left\{\sum_{\ell=0}^{i-j}\left(\prod_{s=1}^\ell\frac{-s+i-j+1}{m-s+i-j+1}\right)d_\ell^{(j)}\right\}\frac{c_j}{\Gamma(i-j+1)}\\
&\quad +\frac{\Gamma(\beta_n+i+m)}{\Gamma(m+1)}\left\{\sum_{\ell=0}^k\left(\prod_{s=1}^\ell\frac{-s+\beta_n+i}{m-s+\beta_n+i}\right)d_\ell^{(n)}\right\}\frac{c_n}{\Gamma(\beta_n+i)}+O(m^{-\alpha_--k-2})\\
&=\sum_{j=1}^i\frac{\Gamma(m+i-j+1)}{\Gamma(m+1)}\left(\prod_{s=1}^{i-j}\frac{-s+i-j+1}{m-s+i-j+1}\right)d_{i-j}^{(j)}\frac{c_j}{\Gamma(i-j+1)}
\\&\quad +\sum_{j=1}^{i}\frac{\Gamma(m+i-j+1)}{\Gamma(m+1)}\left\{\sum_{\ell=0}^{i-j-1}\left(\prod_{s=1}^\ell\frac{-s+i-j+1}{m-s+i-j+1}\right)d_\ell^{(j)}\right\}\frac{c_j}{\Gamma(i-j+1)}\\
&\quad +\frac{\Gamma(\beta_n+i+m)}{\Gamma(m+1)}\left\{\sum_{\ell=0}^k\left(\prod_{s=1}^\ell\frac{-s+\beta_n+i}{m-s+\beta_n+i}\right)d_\ell^{(n)}\right\}\frac{c_n}{\Gamma(\beta_n+i)}+O(m^{-\alpha_--k-2}).
\end{aligned}
\end{equation}
Here it holds
\[
\alpha_-=\min\{-i,-i+1,\ldots,-1,-\re\beta_n-i\}=-i
\]
from the assumptions $n\ge2$ and $\re\beta_n<-n+2$. Now we retake $k\in\mathbb{N}$ sufficiently large as which satisfies $k\ge-\re\beta_n-1$ if needed. Then it holds that
\[
-\alpha_--k-2=i-k-2\le \re\beta_n+i-1.
\]
From this we have
\[
O(m^{-\alpha_--k-2})=O(m^{i-k-2})=O(m^{\re\beta_n+i-1}).
\]
In addition, from Stirling's formula, we have
\[
\frac{\Gamma(\beta_n+i+m)}{\Gamma(m+1)}\left\{\sum_{\ell=0}^k\left(\prod_{s=1}^\ell\frac{-s+\beta_n+i}{m-s+\beta_n+i}\right)d_\ell^{(n)}\right\}\frac{c_n}{\Gamma(\beta_n+i)}=O(m^{\re\beta_n+i-2}).
\]
Let us consider the other terms. On the first term of the right-hand side of \eqref{ind1},
\begin{align*}
\frac{\Gamma(m+i-j+1)}{\Gamma(m+1)}&\left(\prod_{s=1}^{i-j}\frac{-s+i-j+1}{m-s+i-j+1}\right)\frac{1}{\Gamma(i-j+1)}\\
&=\frac{\Gamma(m+i-j+1)}{\Gamma(m+1)}\left(\frac{(i-j)(i-j-1)\cdots\cdot1}{(m+i-j)(m+i-j-1)\cdots\cdot (m+1)}\right)\frac{1}{\Gamma(i-j+1)}\\
&=1
\end{align*}
holds from the property of gamma function.

In the second term of the right-hand side of \eqref{ind1}, by changing the order of summation, we have 
\begin{align*}
\sum_{j=1}^{i}&\frac{\Gamma(m+i-j+1)}{\Gamma(m+1)}\left\{\sum_{\ell=0}^{i-j-1}\left(\prod_{s=1}^\ell\frac{-s+i-j+1}{m-s+i-j+1}\right)d_\ell^{(j)}\right\}\frac{c_j}{\Gamma(i-j+1)}\\
&=\frac{1}{\Gamma(m+1)}\sum_{p=1}^{i-1}\frac{\Gamma(m+p+1)}{\Gamma(p+1)}\sum_{j=1}^{i-p}d_{i-p-j}^{(j)}c_j.
\end{align*}
Therefore we have
\begin{equation}\label{ind2}
\begin{aligned}
\frac{1}{(i-1)!}\sum_{h=0}^m(m-h+1)_{i-1}a_h&=\sum_{j=1}^id_{i-j}^{(j)}c_j+\frac{1}{\Gamma(m+1)}\sum_{p=1}^{i-1}\frac{\Gamma(m+p+1)}{\Gamma(p+1)}\sum_{j=1}^{i-p}d_{i-p-j}^{(j)}c_j\\
&\quad+O(m^{\re\beta_n+i-1}).
\end{aligned}
\end{equation}

Next we consider the left-hand side of \eqref{ind2}. We give a lemma.
\begin{lemma}\label{lem}
The equality
\[
(m-h+1)_\ell=\sum_{p=0}^\ell(-1)^{\ell-p}\binom{\ell}{p}[h]_{\ell-p}(m+1)_p
\]
holds for $\ell\in\mathbb{N}$.
\end{lemma}
This lemma can be proved by induction. From this lemma, we have
\begin{align*}
\frac{1}{(i-1)!}&\sum_{h=0}^m(m-h+1)_{i-1}a_h\\
&=\frac{1}{(i-1)!}\sum_{h=0}^m\left(\sum_{p=0}^{i-1}(-1)^{i-1-p}\binom{i-1}{p}[h]_{i-1-p}(m+1)_p\right)a_h\\
&=\frac{1}{(i-1)!}\sum_{p=0}^{i-1}(-1)^{i-1-p}\binom{i-1}{p}(m+1)_p\sum_{h=0}^m[h]_{i-1-p}a_h\\
&=\frac{(-1)^{i-1}}{(i-1)!}\sum_{h=0}^m[h]_{i-1}a_h+\frac{1}{(i-1)!}\sum_{p=1}^{i-1}(-1)^{i-1-p}\binom{i-1}{p}(m+1)_p\sum_{h=0}^m[h]_{i-1-p}a_h\\
&=\frac{(-1)^{i-1}}{(i-1)!}\sum_{h=0}^m[h]_{i-1}a_h+\frac{1}{\Gamma(m+1)}\sum_{p=1}^{i-1}(-1)^{i-1-p}\frac{\Gamma(m+p+1)}{\Gamma(p+1)\Gamma(i-p)}\sum_{h=0}^m[h]_{i-1-p}a_h.
\end{align*}
Here we used the assumption of induction in the last term and Stirling's formula. Then we obtain
\begin{align*}
\frac{1}{(i-1)!}&\sum_{h=0}^m(m-h+1)_{i-1}a_h\\
&=\frac{(-1)^{i-1}}{(i-1)!}\sum_{h=0}^m[h]_{i-1}a_h+\frac{1}{\Gamma(m+1)}\sum_{p=1}^{i-1}\frac{\Gamma(m+p+1)}{\Gamma(p+1)}\left(\sum_{j=1}^{i-p}d_{i-p-j}^{(j)}c_j+O(m^{\re\beta_n+i-p-1})\right)\\
&=\frac{(-1)^{i-1}}{(i-1)!}\sum_{h=0}^m[h]_{i-1}a_h+\frac{1}{\Gamma(m+1)}\sum_{p=1}^{i-1}\frac{\Gamma(m+p+1)}{\Gamma(p+1)}\sum_{j=1}^{i-p}d_{i-p-j}^{(j)}c_j+\sum_{p=1}^{i-1}\frac{1}{\Gamma(p+1)}O(m^{\re\beta_n+i-1})\\
&=\frac{(-1)^{i-1}}{(i-1)!}\sum_{h=0}^m[h]_{i-1}a_h+\frac{1}{\Gamma(m+1)}\sum_{p=1}^{i-1}\frac{\Gamma(m+p+1)}{\Gamma(p+1)}\sum_{j=1}^{i-p}d_{i-p-j}^{(j)}c_j+O(m^{\re\beta_n+i-1}).
\end{align*}
Substituting this into \eqref{ind2}, we have
\[
\begin{aligned}
\frac{(-1)^{i-1}}{(i-1)!}\sum_{h=0}^m[h]_{i-1}a_h&=\sum_{j=1}^id_{i-j}^{(j)}c_j+O(m^{\re\beta_n+i-1}).
\end{aligned}
\]
This completes the proof of Proposition \ref{c1-cn-1}.

\subsection{Proof of Proposition \ref{hgs_conv}}
In the following, we give a proof of Proposition \ref{hgs_conv}. From the definition of $[h]_{i-1}$, we obtain
\begin{align*}
\sum_{h=0}^m[h]_{i-1}a_h&=\sum_{h=0}^mh(h-1)\cdots(h-i+2)\frac{(\bm{\alpha}_0)_h}{(\bm{\beta}_0)_hh!}\\
&=\sum_{h=i-1}^m\frac{(\bm{\alpha}_0)_h}{(\bm{\beta}_0)_h(h-i+1)!}
\end{align*}
for $1\le i\le n-1$. We set $\ell=h-i+1$. Then we have
\begin{align*}
\sum_{h=i-1}^m\frac{(\bm{\alpha}_0)_h}{(\bm{\beta}_0)_h(h-i+1)!}&=\sum_{\ell=0}^{m-i+1}\frac{(\bm{\alpha}_0)_{\ell+i-1}}{(\bm{\beta}_0)_{\ell+i-1}\ell!}.
\end{align*}
Since
\begin{align*}
(a)_{\ell+i-1}&=a(a+1)\cdots(a+i-2)(a+i-1)\cdots(a+i+\ell-2)\\
&=(a)_{i-1}(a+i-1)_\ell
\end{align*}
holds for $a\in\mathbb{C}$, we obtain
\[
\frac{(\bm{\alpha}_0)_{\ell+i-1}}{(\bm{\beta}_0)_{\ell+i-1}}=\frac{(\bm{\alpha}_0)_{i-1}}{(\bm{\beta}_0)_{i-1}}\frac{(\bm{\alpha}_0+i-1)_{\ell}}{(\bm{\beta}_0+i-1)_{\ell}}.
\]
Then we have
\[
\sum_{h=0}^m[h]_{i-1}a_h=\frac{(\bm{\alpha}_0)_{i-1}}{(\bm{\beta}_0)_{i-1}}\sum_{\ell=0}^{m-i+1}\frac{(\bm{\alpha}_0+i-1)_{\ell}}{(\bm{\beta}_0+i-1)_{\ell}\ell!}.
\]
From Raabe's ratio test, we see that the power series
\[
\sum_{\ell=0}^{\infty}\frac{(\bm{\alpha}_0+i-1)_{\ell}}{(\bm{\beta}_0+i-1)_{\ell}\ell!}
\]
is convergent when $\re\beta_n<-n+2$, and its value is nothing but ${}_nF_{n-1}(\bm{\alpha}_0+i-1,\bm{\beta}_0+i-1;1)$. In conclusion, we have
\[
\lim_{m\to\infty}\sum_{h=0}^m[h]_{i-1}a_h=\frac{(\bm{\alpha}_0)_{i-1}}{(\bm{\beta}_0)_{i-1}}{}_nF_{n-1}(\bm{\alpha}_0+i-1,\bm{\beta}_0+i-1;1).
\]

\section*{Acknowledgement}
The author would like to thank Professor Yoshishige Haraoka for his valuable comments and suggestions.

\medskip

\begin{flushleft}
Shunya Adachi\\
Graduate School of Science and Technology\\
Kumamoto University\\
2-39-1 Kurokami, Chuo-ku, Kumamoto, 860-8555, Japan\\
E-mail address: {\tt 200d7101@st.kumamoto-u.ac.jp}
\end{flushleft}

\end{document}